\documentclass[12pt]{amsart}
\usepackage{amscd}
\def\frk{\frak}               

\def\mm{{\frk m}}

\def\Phi{{\frk n}}
\def\Phi{{\frk N}}
\def\opn#1#2{\def#1{\operatorname{#2}}} 
\opn\chara{char} \opn\length{\ell} \opn\pd{pd} \opn\rk{rk}
\opn\projdim{proj\,dim} \opn\injdim{inj\,dim} \opn\rank{rank}
\opn\depth{depth} \opn\grade{grade} \opn\height{height}
\opn\embdim{emb\,dim} \opn\codim{codim}

\opn\Tr{Tr} \opn\bigrank{big\,rank}
\opn\superheight{superheight}\opn\lcm{lcm}
\opn\trdeg{tr\,deg}%
\opn\reg{reg} \opn\lreg{lreg} \opn\ini{in} \opn\lpd{lpd}
\opn\size{size}
\opn\div{div} \opn\Div{Div} \opn\cl{cl} \opn\Cl{Cl}
\opn\Spec{Spec} \opn\Supp{Supp} \opn\supp{supp} \opn\Sing{Sing}
\opn\Ass{Ass}
\opn\Ann{Ann} \opn\Rad{Rad} \opn\Soc{Soc}
\opn\Im{Im} \opn\Ker{Ker} \opn\Coker{Coker} \opn\Am{Am}
\opn\Hom{Hom} \opn\Tor{Tor} \opn\Ext{Ext} \opn\End{End}
\opn\Aut{Aut} \opn\id{id}

\opn\nat{nat}
\opn\pff{pf}
\opn\Pf{Pf} \opn\GL{GL} \opn\SL{SL} \opn\mod{mod} \opn\ord{ord}
\opn\Gin{Gin} \opn\Hilb{Hilb}
\opn\aff{aff} \opn\con{conv} \opn\relint{relint} \opn\st{st}
\opn\lk{lk} \opn\cn{cn} \opn\core{core} \opn\vol{vol}
\opn\link{link} \opn\star{star}
\opn\gr{gr}

\def\pot#1#2{#1[\kern-0.28ex[#2]\kern-0.28ex]}

\opn\dirlim{\underrightarrow{\lim}}
\opn\inivlim{\underleftarrow{\lim}}
\let\tensor=\otimes

\let\Dirsum=\bigoplus

\let\to=\rightarrow
\let\To=\longrightarrow
\def\Implies{\ifmmode\Longrightarrow \else
        \unskip${}\Longrightarrow{}$\ignorespaces\fi}
\def\implies{\ifmmode\Rightarrow \else
        \unskip${}\Rightarrow{}$\ignorespaces\fi}
\def\iff{\ifmmode\Longleftrightarrow \else
        \unskip${}\Longleftrightarrow{}$\ignorespaces\fi}

\let\:=\colon
\newtheorem{Theorem}{Theorem}[section]
\newtheorem{Lemma}[Theorem]{Lemma}

\newtheorem{Remarks}[Theorem]{Remarks}

\let\epsilon\varepsilon
\let\phi=\varphi
\let\kappa=\varkappa
\def\qed{\ifhmode\textqed\fi
      \ifmmode\ifinner\quad\qedsymbol\else\dispqed\fi\fi}
\def\textqed{\unskip\nobreak\penalty50
       \hskip2em\hbox{}\nobreak\hfil\qedsymbol
       \parfillskip=0pt \finalhyphendemerits=0}
\def\dispqed{\rlap{\qquad\qedsymbol}}

\opn\dis{dis}
\def\pnt{{\raise0.5mm\hbox{\large\bf.}}}
\def\lpnt{{\hbox{\large\bf.}}}
\opn\Lex{Lex}

\begin{document}

\title{Notes  on the multiplicity conjecture}

\author{J\"urgen Herzog and  Xinxian Zheng}
\address{J\"urgen Herzog, Fachbereich Mathematik und
Informatik, Universit\"at Duisburg-Essen, Campus Essen, 45117
Essen, Germany} \email{juergen.herzog@uni-essen.de}

\address{Xinxian Zheng, Fachbereich Mathematik und
Informatik, Universit\"at Duisburg-Essen, Campus Essen, 45117
Essen, Germany} \email{xinxian.zheng@uni-essen.de} \maketitle

\begin{abstract}
New cases of the multiplicity conjecture are considered.
\end{abstract}

\section*{Introduction}
Throughout this paper, we fix a field $K$ and let $R$ be a
homogeneous $K$-algebra. In other words, $R$ is a finitely
generated $K$-algebra, generated over $K$ by elements of degree
$1$, and hence is isomorphic to $S/I$ where $S=K[x_1,\ldots, x_n]$
is a polynomial ring and $I$ a graded ideal contained in
$(x_1,\ldots, x_n)$. Consider a graded minimal free $S$-resolution
of $R$:
\[
0\To\bigoplus_{j=1}^{b_p} S(-d_{pj})\To\cdots\To
\bigoplus_{j=1}^{b_1} S(-d_{1j})\To S \To 0.
\]
The ring $R$ is said to have a {\em pure resolution} if for all
$i$, the shifts $d_{ij}$ do not depend on $j$ (but only on $i$).
Hence if the resolution is pure, it has the following shape:
\[
0\To S^{b_p}(-d_{p})\To\cdots\To S^{b_1}(-d_{1})\To S \To 0.
\]
When $R$ is Cohen-Macaulay and has a pure resolution, Huneke and
Miller's formular \cite{HM} says that the multiplicity  of $R$ is
given by
\[
e(R)=(\prod_{i=1}^p d_i)/p!.
\]
In general we define $M_i(I)=\max\{d_{ij}\, : j=1,\ldots, b_i\}$
and $m_i(I)=\min\{d_{ij}\, : j=1,\ldots, b_i\}$ for $i=1,\ldots
p$. When there is no danger of ambiguity, we write $M_i$ and $m_i$
for short.

Huneke and Srinivasan had the following

\medskip\noindent
{\bf Conjecture 1.}\ {\em  For each homogeneous Cohen-Macaulay
$K$-algebra $R$
\[
(\prod_{i=1}^p m_i)/p!\leq e(R) \leq (\prod_{i=1}^p M_i)/p!.
\]
}

\medskip
\noindent
 Conjecture 1 has been widely studied and partial results
have been obtained. In \cite{HS}, the first author and Srinivasan
showed that this conjecture is true in the following cases: $R$
has a quasi-pure resolution (i.e., $m_i(I)\geq M_{i-1}(I)$); $I$
is a perfect ideal of codimension $2$; $I$ is a codimension $3$
Gorenstein ideal generated by $5$ elements; $I$ is a stable ideal;
$I$ is a squarefree strongly stable ideal. Furthermore, Guardo and
Van Tuly \cite{GT} proved that the conjecture holds for powers of
complete intersections, Srinivasan \cite{S} proved a stronger
bounds for Gorenstein algebras with quasi-pure resolutions. And
recently, in \cite{MNR}, Migliore, Nagel and R\"omer proved  a
stronger version of Conjecture  1 when $R$ is a codimension $2$ or
Gorenstein codimension $3$ algebra (with no limitations on the
number of generators). As a corollary, they showed that in these
two cases, the multiplicity $e(R)$ reaches the upper or lower
bound if and only if $R$ has a pure resolution.

It was observed in \cite{HS}, that the lower bound in Conjecture 1
fails in general if $R$ is not Cohen-Macaulay, even one replaces
Conjecture 1 the projective dimension by the codimension. In the
same paper, the authors had the following stronger conjecture for
the upper bound of the multiplicity of $R$.

\medskip
\noindent {\bf Conjecture 2.}\ {\em  Let $R$ be a homogeneous
$K$-algebra of codimension $s$. Then
\[
e(R)\leq (\prod_{i=1}^s M_i)/s!.
\]}

\medskip
\noindent If the defining ideal of $R$ is stable, or squarefree
strongly stable, or if $R$ has a linear resolution, Conjecture 2
is shown to be true in \cite{HS}. In addition, Gold \cite{G}
proved it for codimension $2$ lattice ideals. This was generalized
by R\"omer \cite{R} for all codimension two ideals. In the same
paper  R\"omer proved Conjecture 2 for componentwise linear ideals
in characteristic 0.

In Section 1 we show that if $I\subset S$ is an ideal of
codimension $s$, not necessarily perfect, for which one has
\[
(\prod_{i=1}^s m_i)/s!\leq e(R) \leq (\prod_{i=1}^s M_i)/s!,
\]
and if $f_1,\ldots, f_m$ is a regular sequence modulo $I$, then
the corresponding inequalities are again valid for $(I,
f_1,\ldots, f_m)$.  One might expect that the proof of this
statement is rather simple. But again careful estimates are
required to establish the result.

In  Section 2 we show that Conjecture  2 is valid in the limit
with respect to taking powers of ideal, that is, we show that
\[
\lim_{k\to
\infty}\frac{e(S/I^k)}{\frac{1}{s!}\prod_{i=1}^sM_i(I^k)}\leq 1.
\]
Unfortunately, this does not imply that Conjecture  2 holds for
all sufficiently high powers of $I$.

In view of the  results by Migliore, Nagel and R\"omer \cite{MNR}
one is lead to ask the following question: suppose that for a ring
$R$ the lower bound in Conjecture  1, or the upper bound in
Conjecture 2 is reached. Does this imply that $R$ is
Cohen-Macaulay and  has a pure resolution?

In Section 3 we show that  not only in the cases described by
Migliore, Nagel and R\"omer this improved multiplicity conjecture
holds, but also for rings with almost pure resolutions, as well as
for rings defined by componentwise linear ideals.

We also recall a recent result of Mir\'o-Roig \cite{M} who showed
that graded ideals of maximal minors of maximal grade  satisfy
Conjecture 1, and we show that again in this case the upper or
lower bound is reached if and only if the resolution is pure. This
case includes rings whose defining ideal is a power of a graded
regular sequence.

\section{The multiplicity conjecture and regular sequences}
Let $I\subset S$ be a graded ideal of codimension $s$, $R=S/I$ and
$f_1, \ldots f_m$ a homogeneous regular sequence of $R$. Suppose
that \[ (\prod_{i=1}^s m_i(I))/s!\leq e(R) \leq (\prod_{i=1}^s
M_i(I))/s!. \]

Are the corresponding inequalities again valid for  for $R/(f_1,
\ldots, f_m)$?  In fact,   without assuming that $S/I$ is
Cohen-Macaulay, we have the following
\begin{Theorem}
\label{regular} Let $I\subset S$ be a graded ideal of codimension
$s$, $R=S/I$ and ${\bf f}=f_1, \ldots f_m$ a homogeneous regular
sequence of $R$. Suppose that
\[
(\prod_{i=1}^s m_i(I))/s!\leq e(R) \leq (\prod_{i=1}^s M_i(I))/s!.
\] Then
\[
(\prod_{i=1}^{s+m} m_i(I, {\bf f}))/(s+m)!\leq e(R/({\bf f})) \leq
(\prod_{i=1}^{s+m} M_i(I,{\bf f}))/(s+m)!.
\]
\end{Theorem}

\begin{proof}
By using induction on $m$, one needs only to show the case $m=1$.
For simplicity, we denote $f_1$ by $f$, $M_i(I)$ and $m_i(I)$ by
$M_i$ and $m_i$, respectively. Let $d=\deg f$.

We first show that $e(R/(f)) \leq (\prod_{i=1}^{s+1}
M_i(I,f))/(s+1)!$. We have $\codim (I,f)=s+1$ and
$e(R/(f))=e(S/(I,f))=e(R)\cdot d.$

 Let $G_\lpnt$ be the minimal graded free resolution of
$R$, and $H_\lpnt$ the minimal graded free resolution of $S/(f)$.
Then $F_\lpnt=G_\lpnt\otimes H_\lpnt$ is the minimal graded free
resolution of $R/(f)=S/(I,f)$. Hence $F_i=G_i\otimes S\oplus
G_{i-1}\otimes R(-d)$, $i=1,\ldots, p+1$, where $p=\projdim S/I$.
Therefore $M_i((I,f))=\max\{M_i, M_{i-1}+d\}$ for $i=1,\ldots,
p+1$.

Thus  we need to show that
\begin{eqnarray}
\label{*} (s+1)d\prod_{i=1}^s  M_i \leq \prod_{i=1}^{s+1} \max
\{M_i, M_{i-1}+d\}, \text{ where } M_0=0.
\end{eqnarray}
Set $M_i=id+y_i$ for $1\leq1\leq s+1$. Then $\max \{M_i,
M_{i-1}+d\}=\max \{id+y_i, id+y_{i-1}\}=id+\max\{y_{i-1}, y_i\}$.
Let $N=\{i\, : y_i\geq 0, \quad 1\leq i\leq s+1\}$, and let
$j=\max\{i\, : i\in N\}$. In case $N=\emptyset$, we set $j=0$.
Then $y_i<0$ for all $i$ with $j<i\leq s+1$.

We will distinguish two cases:

\medskip
\noindent Case 1. $j=s+1$. We have $\max\{y_s, y_{s+1}\}\geq
y_{s+1}\geq 0$. Inequality (\ref{*}) is equivalent to the
inequality
\[
(s+1)d\prod_{i=1}^s M_i \leq \prod_{i=1}^s \max \{M_i,
M_{i-1}+d\}\cdot((s+1)d+\max\{y_s, y_{s+1}\}),
\]
which is satisfied since $M_i\leq \max\{M_i, M_{i-1}+d\}$ for
$i=1,\ldots, s$, and because
 $\max\{y_s, y_{s+1}\}\geq 0$.

\medskip
\noindent
 Case 2. $j<s+1$. It is suffices to show that
\begin{eqnarray}
\label{**} \;\; \;\;\; \prod_{i=j+1}^s (id+y_i)\cdot (s+1)d\leq
\prod_{i=j+1}^{s+1}
(id+y_{i-1})=((j+1)d+y_j)\prod_{i=j+2}^{s+1}(id+y_{i-1}).
\end{eqnarray}
Set $z_i=y_i/d$. Then inequality (\ref{**}) is equivalent to the
inequality
\[
\prod_{i=j+1}^s(i+z_i)\cdot(s+1)
\leq((j+1)+z_j)\prod_{i=j+1}^s((i+1)+z_i).
\]
Since $id+y_i=M_i>0$, it follows $0<i+z_i$ for $i=1,\ldots,s+1$.
Hence we need to show that
\[
s+1\leq ((j+1)+z_j)\cdot \prod_{i=j+1}^s(1+1/{(i+z_i)}).
\]
Since $z_j=y_j/d\geq 0$ and $z_i<0$ for all $i=j+1,\ldots,s$, and
since $i+z_i=i+y_i/d=(id+y_i)/d=m_i/d>0$, it follows that
$1/{(i+z_i)}>1/i$. Therefore
\[
s+1=(j+1)\prod_{i=j+1}^s(1+1/i)\leq((j+1)+z_j)\prod_{i=j+1}^s(1+1/(i+z_i)).
\]

Similarly, by taking $j=\max\{i\,: y_i\leq 0, \quad 0\leq i\leq
s+1\}$, where $y_0=0$, and distinguishing the cases $j=s+1$ and
$j<s+1$, one sees that $e(R/(f)) \geq (\prod_{i=1}^{s+1}
m_i(I,f))/(s+1)!$.
\end{proof}

\section{Powers  of an ideal}

As the main result of this section we want to prove that
Conjecture  2 is true in the limit with respect to taking powers
of an ideal. To be more precise, we show

\begin{Theorem}
\label{limit} For any graded ideal $I\subset S$ of codimension $s$
we have
\[
\lim_{k\to
\infty}\frac{e(S/I^k)}{\frac{1}{s!}\prod_{i=1}^sM_i(I^k)}\leq 1.
\]
\end{Theorem}

\begin{proof}
For the proof of this result we will proceed in several steps.

\medskip
 (i) Let $M$ be a graded $S$-module of projective dimension $p$. We set
\[
\reg_i(M)=\max\{j\: \beta_{ii+j}(M)\neq 0\}\quad \text{for}\quad
i=0,\ldots,p.
\]
Then $M_{i}(I)=\reg_i(S/I)+i$ for  $i=1,\ldots, p$, and
$$\reg(I)=\max\{\reg_i(I)\: i=0,\ldots, p\}$$ is the regularity of
$I$.

Let $L\subset S$ be a graded ideal of codimension $s$. We claim
that  $\reg_i(S/L)\geq \reg_{i-1}(S/L)$ for $i=0,\ldots, s$. In
fact, let
\[
\begin{CD}
0@>>> F_p@>\phi_p >> F_{p-1}@>\phi_{p-1} >> \cdots @>\phi_2 >>
F_1@>\phi_1 >> F_0@>>> S/L@>>> 0
\end{CD}
\]
be a graded minimal free resolution of $S/L$. Suppose that
$\reg_i(S/L)<\reg_{i-1}(S/L)$ for some $i\leq s$. Then
$M_{i}(L)\leq M_{i-1}(L)$. Let $e\in F_{i-1}$ be a homogeneous
basis element of degree $M_{i-1}(L)$, and let $f$ be any
homogeneous basis element of $F_i$. Then $\deg \phi_i(f)=\deg
f\leq M_{i}(L)\leq M_{i-1}(L)$. Thus if we write  $\phi_i(f)$ as a
linear combination of the basis elements of $F_{i-1}$, the
coefficient $a$  of $e$ will be of degree $\deg a\leq
M_{i}(L)-M_{i-1}(L)\leq 0$. This is only possible if $a=0$, since
$\phi_i(f)\in \mm F_{i-1}$, where $\mm$ is the graded maximal
ideal of $S$.

Now we consider the $S$-dual of the resolution $F$. Let $e^*$ be
the dual basis element of $e$, and $\phi_i^*\: F_{i-1}^*\to F_i^*$
the map dual to $\phi_i$. Then it follows that $\phi_i^*(e^*)=0$.
Thus $e^*$ is a cycle of the dual complex. On the other hand,
$e^*$ cannot be a boundary, since $e^*$ is a basis element of
$F_{i-1}^*$ and since the image of $\phi_{i-1}^*\: F_{i-2}^*\to
F_{i-1}^*$ is contained in $\mm F_{i-1}^*$. Hence we see that
$\Ext_S^{i-1}(S/L,S)=H^{i-1}(F^*)\neq 0$. This contradicts the
fact that  $\Ext_S^j(S/L,S)=0$ for $j<s$, since $\grade L=\codim
L=s$.

\medskip (ii) Cutkosky, Herzog and Trung \cite{CHT} as well as
Kodiyalam \cite{K} showed that $$\reg_i(I^k)=q_ik+c_i$$ is a
linear function for $k\gg 0$. In particular, $\reg(I^k)=qk+c$ for
$k\gg 0$. It is also shown that $q=q_0$, see \cite[Corollary
3.2]{CHT}.

By (i) we have
\[
\reg_0(I^k)\leq \reg_i(I^k)\leq \reg(I^k)
\]
for $i=0,\ldots,s-1$ where $s=\codim I$ ($=\codim I^k$ for all
$k$). Thus (2) implies
\[
qk+c_0\leq q_ik+c_i\leq qk+c
\]
for $i=0,\ldots,s-1$ and all $k\gg 0$. This implies that $q_i=q$
for $i=0,\ldots,s-1$.

\medskip
(iii) From (ii) it  follows that
\[
\frac{1}{s!}\prod_{i=1}^s M_i(I^k)=\frac{q^s}{s!}k^s+\cdots \quad
\text{for}\quad k\gg 0.
\]
is a polynomial function of degree $s$ whose leading coefficient
is $q^s/s!$.
\medskip

(iv) For $k\gg 0$, the function $s!e(S/I^k)$ is a polynomial
function whose leading term is an integer which we denote by
$e(I,S)$. In other words,
\[
e(S/I^k)=\frac{e(I,S)}{s!}k^s+\cdots  \quad \text{for}\quad k\gg
0.
\]
 For all $k\gg 0$, let $\{P_1,\ldots, P_r\}$ be the
(stable) set of minimal prime ideals of $S/I^k$. The associativity
formula of multiplicity (\cite[Corollary 4.7.8]{BH}) then shows
that
\[
e(S/I^k)=\sum_{i=1}^r\length(S_{P_i}/I_{P_i}^k)e(S/P_i).
\]
Here $\length(M)$ denotes the length of a module $M$.

Each $S_{P_i}$ is a regular local ring of dimension $h$, and
$I_{P_i}$ is $P_iS_{P_i}$ primary. Therefore
\begin{eqnarray*}
\length(S_{P_i}/I_{P_i}^k)=\frac{e(I_{P_i},S_{P_i})}{s!}k^s+\cdots
\end{eqnarray*}
is a polynomial function of degree $s$ for $k\gg 0$, see
\cite[Proposition 4.2.6]{BH}, where the numerator
$e(I_{P_i},S_{P_i})$ of the leading coefficient of this polynomial
is  the multiplicity of $S_{P_i}$ with respect to $I_{P_i}$. This
proves our assertion and also shows that
\begin{eqnarray}
\label{formula}
 e(I,S)=\sum_Pe(I_P,S_P)e(S/P),
\end{eqnarray}
where the sum is taken over all asymptotic minimal prime ideals of
$I$.

\medskip
\noindent (v) The theorem will follow once we have shown
$e(I,S)\leq q^s$. We first notice that $e(I,S)\leq e(J, S)$ for
any
 ideal $J\subset I$ with $\codim J=\codim I=s$,  and that $e(I,S)= e(J, S)$ if $J$ is a reduction ideal of $I$, that is,
 if $JI^k=I^{k+1}$ for some $k$. Indeed, this follows from formula
(\ref{formula}) and the fact that the corresponding statements are
true for ideals in a local ring which are primary to the maximal
ideal, see \cite[Lemma 4.6.5]{BH}.

Now we use the fact, shown by Kodiyalam \cite[Theorem 5]{K}, that
$I$ admits a reduction ideal $J$ with $\reg_0(J)=q$. Hence
replacing $I$ by $J$ we may assume that $I$ is generated in degree
$\leq q$.

After a base field extension we may assume that $K$ is infinite.
Then, since $\codim I=s$,  generically chosen $q$-forms
$f_1,\ldots, f_s\in I_q$ will form a regular sequence. Let $L$ be
the ideal generated by these forms. Then $e(I,S)\leq e(L,S)=q^s$,
as desired. It just remains to establish the last equation. This
can be seen as follows: Since the $L$ is generated by a regular
sequence each of the factor modules $L^{k-1}/L^{k}$ is a free
$S/L$-module, whose rank is $\binom{s+k-1}{s-1}$. Hence, since
$e(S/L)=q^s$, we see that
$e(L^{k-1}/L^{k})=q^s\binom{s+k-1}{s-1}$. It follows that
\[
e(S/L^k)=q^s\sum_{j=1}^k\binom{s+k-1}{s-1}=\frac{q^s}{s!}k^s+\cdots
\]
This implies that $e(L,S)=q^s$.
\end{proof}

Unfortunately Theorem \ref{limit} does not imply that Conjecture 2
is true for all high enough powers of an ideal, as it is easy to
find  ideals for which $$\lim_{k\to
\infty}\frac{e(S/I^k)}{\frac{1}{s!}\prod_{i=1}^sM_i(I^k)}=1.$$

\section{The improved multiplicity conjecture}

Motivated by the results of Migliore, Nagel and R\"omer
\cite{MNR}, we say that the improved multiplicity conjecture
holds, if all standard graded $K$-algebras $R$ satisfy the
multiplicity conjectures, and whenever the bounds are reached,
then the defining ideal has a pure resolution and $R$ is
Cohen-Macaulay.

In this section we show  that for some interesting classes of
examples the improved multiplicity conjecture holds.

Generalizing the result \cite[Corollary 1.3]{MNR} we first show

\begin{Theorem}
\label{roemer} Let $I\subset S$ be a graded ideal of codimension
$2$. Then $S/I$ satisfies the improved multiplicity conjecture.
\end{Theorem}

\begin{proof}
R\"omer proved  in \cite[Theorem 2.4]{R}  that $R=S/I$ satisfies
Conjecture 2. Thus it remains to be shown that if
$e(R)=(1/2)M_1M_2$, then $R$ is Cohen-Macaulay and has a pure
resolution. Once it is shown that $R$ is Cohen-Macaulay, then by
\cite[Theorem ]{MNR} we also have that  $R$ has pure resolution.
(This last fact also follows from Theorem \ref{holds} below.)

Let $S=K[x_1,\ldots,x_n]$. One may assume that $K=\infty$, and
that (after a generic change of coordinates)  $x_1,\ldots, x_n$ is
an almost regular sequence on $R$, i.e, multiplication with $x_i$
on $R_{i-1}=R/(x_1,\ldots,x_{i-1})R$ has a finite length kernel
for all $i$. In his proof,  R\"omer showed that
\[
e(R)\leq e(R_{n-2})\leq (1/2)M_1M_2,.
\]
He also showed in \cite[Lemma 2.3]{R} that $e(R_{n-2})
=e(R)+\text{length}(0:_{R_{n-3}}x_{n-2})$. Thus if we assume that
the upper bound is reached, then $x_{n-2}$ is regular on
$R_{n-3}$. By \cite[Proposition 3]{HK}  this implies that
$x_1,\ldots, x_{n-2}$ is a regular sequence, and hence $R$ is
Cohen-Macaulay.
\end{proof}

We say that $I$ is {\em componentwise linear} if all ideals
spanned by the graded components of $I$ have a linear resolution.
Our next  class of rings satisfying the improved multiplicity
conjecture is the following

\begin{Theorem}
\label{componentwise} Let $I\subset S$ be a componentwise linear
ideal of codimension $s$. Then for $S/I$ the improved multiplicity
conjecture holds.
\end{Theorem}

Before we prove this theorem we first note

\begin{Lemma}
\label{purelinear} Let $I\subset S$ be a componentwise linear
ideal with a pure resolution. Then $I$ has a linear resolution.
\end{Lemma}

\begin{proof}
We may assume that the base field is infinite. Let $\Gin(I)$
denote the generic initial ideal of $I$ with respect to the
reverse lexicographical order. In \cite{AHH} it is shown that $I$
is componentwise linear if and only if $I$ and $\Gin(I)$ have the
same graded Betti-numbers, provided $\chara K=0$, and that
$\Gin(I)$ is a strongly stable ideal. Here we only need that $I$
and $\Gin(I)$ have the same graded Betti-numbers and that
$\Gin(I)$ is stable. The proof given in \cite{AHH} shows that this
is the case in all characteristics. Indeed, let $I_{\langle
j\rangle}$ be the ideal generated by all elements of degree $j$ in
$I$. By our assumption on $I$, the ideal $I_{\langle j\rangle}$
has a linear resolution. Applying the Bayer-Stillman theorem
\cite{BS} it follows that $\Gin(I_{\langle j\rangle})$ has a
linear resolution, so that $\Gin(I_{\langle
j\rangle})=\Gin(I)_{\langle j\rangle}$. This proves that $\Gin(I)$
is again componentwise linear. Since $\Gin(I_{\langle
j\rangle})=\Gin(I)_{\langle j\rangle}$ is $p$-Borel and has a
linear resolution, Proposition  10 of Eisenbud, Reeves and Totaro
\cite{ERT} implies that $\Gin(I)_{\langle j\rangle}$ is a  stable
 monomial ideal for all $j$, and hence $\Gin(I)$ is a stable monomial
 ideal.
 With the same arguments as in the proof of \cite[Theorem 1.1]{HH}
 it then follows that $I$ and $\Gin(I)$ have the same graded
 Betti-numbers.

Replacing $I$ by $\Gin(I)$ we may as well assume that $I$ is a
stable ideal. The Eliahou-Kervaire resolution \cite{EK} of $S/I$
implies that
\[
M_i=\max\{\deg u\: u\in G(I),\quad m(u)\geq i\}+i-1\}, \] and
\[
m_i=\min\{\deg u\: u\in G(I),\quad m(u)\geq i\}+i-1\}.
\]
Here $G(I)$ is the unique minimal set of monomial generators of
$I$, and $m(u)$ denote for a monomial $u$ the largest integer $i$
such that $x_i|u$.

Thus if $m_i=M_i$ for all $i$, then $I$ is generated in one degree
and hence has a linear resolution, as it is a stable ideal.
\end{proof}

\begin{proof}[Proof Theorem \ref{componentwise}]
As shown in the preceding  lemma, we may assume that $I$ is a
stable  monomial ideal. Its Betti-numbers do not depend on the
characteristic of the base field. Thus we may assume that the base
field has characteristic $0$. Since $I$ is componentwise linear it
follows from \cite[Theorem 1.1.]{HH} (see also the proof of Lemma
\ref{purelinear}) that $\Gin(I)$ is again componentwise linear and
that $I$ and $\Gin(I)$ have the same graded Betti-numbers.
Replacing $I$ by $\Gin(I)$ and observing that $\Gin(I)$ is
strongly stable since the characteristic of the base field is $0$,
we may now assume that $I$ is strongly stable.

The proof of the multiplicity conjecture for stable ideals given
in \cite[Theorem 3.2]{HS} is in fact only valid for strongly
stable ideals, as it is used there that if $I\subset
K[x_1,\ldots,x_n]$ is stable, then $I:x_n$ is stable as well. But
this is only true for strongly stable ideals. However, as seen
above, the stable ideal may be replaced by a strongly stable
ideal. Thus for the rest of our proof  we may follows the
arguments given the proof of \cite[Theorem 3.2]{HS}.

We first treat the case that $S/I$ is Cohen-Macaulay. In that case
we may assume that $x_n^a\in I$ where $n=\dim_KS_1$, so that in
particular $S/I$ is Artinian. In the proof of \cite[Theorem
3.2]{HS} it is shown that $e(S/I)\leq e(S/(x_1,\ldots,
x_n)^a)=(1/n!)\prod_{i=1}^nM_i$. Thus if the upper bound is
reached, then $I=(x_1,\ldots, x_n)^a$, and so $S/I$ has a pure
resolution.

Now suppose $S/I$ reaches the lower bound. We prove the assertion
by induction on the length of $S/I$. The case
$\text{length}(S/I)=1$ is trivial. So now we assume that
$\text{length}(S/I)>1$. Let $J\subset \bar{S}=K[x_1,\ldots,
x_{n-1}]$ be the unique monomial ideal such that
$(J,x_n)=(I,x_n)$. The ideals $J$ and $(I:x_n)$ are again strongly
stable ideals, and since the multiplicity conjecture holds for
strongly stable ideals we have
\begin{eqnarray*}
e(S/I)&=&e(S/(I,x_n))+e(S/(I:x_n))\\
&\geq &
(1/(n-1)!)\prod_{i=1}^{n-1}m_i(J)+(1/n!)\prod_{i=1}^nm_i(I:x_n).
\end{eqnarray*}
It is shown in the proof of \cite[Theorem 3.2]{HS} that the right
hand side of this inequality  is greater that or equal to
$(1/n!)\prod_{i=1}^nm_i$.

Our assumption implies that
$e(S/(I,x_n))=e(\bar{S}/J)=(1/(n-1)!)\prod_{i=1}^{n-1}m_i(J)$ and
$e(S/(I:x_n))=(1/n!)\prod_{i=1}^nm_i(I:x_n)$. Hence the induction
hypothesis yields that both $\bar{S}/J$ and $S/(I:x_n)$ have a
pure resolutions. Lemma \ref{purelinear} then implies that both
$\bar{S}/J$ and $S/(I:x_n)$ have a linear  resolution. Since $S/I$
is Artinian, it follows that $\bar{S}/J$ and $S/(I:x_n)$ are
Artian, and so there exist numbers $a$ and $b$ such that
$J=(x_1,\ldots, x_{n-1})^a$, and $(I:x_n)=(x_1,\ldots,x_n)^b$.
Therefore, $$I=(x_1,\ldots, x_{n-1})^a+(x_1,\ldots,x_n)^bx_n.$$ If
$n=1$, then $I$ has a linear resolution. Thus we now may assume
that $n>1$.

Since $J\subset I\subset (I:x_n)$ it follows that $a\geq b$, and
since $I$ is strongly stable it follows that $a\leq b+1$. Suppose
that $a=b$. Then
\begin{eqnarray*}
e(S/I)&=&e(\bar{S}/J)+e(S/(I:x_n)\\
&=&{n+a-2\choose n-1}+{n+a-1\choose n},
\end{eqnarray*}
so that $$n!e(S/I)=(2n+a-1)\prod_{i=0}^{n-2}(a+i).$$

On the other hand, $m_i=a+i-1$ for $i=1,\ldots,n-1$ and $m_n=a+n$.
Therefore, $n!e(S/I)\neq \prod_{i=1}^nm_i$ for $n>1$, a
contradiction. Hence we conclude that $a=b+1$, and so
$I=(x_1,\ldots, x_n)^{b+1}$. In particular, $I$ has a linear
resolution.

Finally, let $I$ be an arbitrary strongly stable ideal such that
the multiplicity of $S/I$ reaches the upper bound. We want to show
that $S/I$ is Cohen-Macaulay.

We may assume that $(I:x_n)\neq I$, because otherwise $I\subset
\bar{S}=K[x_1,\ldots, x_{n-1}]$, and we are done by induction on
$n$. Assume $S/I$ is not Cohen-Macaulay. In this case it shown in
the proof \cite[Theorem 3.7]{HS} that
\[
e(S/I)=e(S/(I:x_n))\leq 1/t!\prod_{i=1}^tM_i(I:x_n)\leq
1/s!\prod_{i=1}^sM_i,
\]
where $t=\codim S/(I:x_n)\leq \codim S/I=s$. It is also shown in
\cite[Lemma 3.6]{HS} that $M_i(I:x_n)\leq M_i$ for all $i$. Since
on the other hand, $M_i(I)/i\geq 1$ for all $i$,   our assumption
implies that
\begin{enumerate}
\item[(i)] $t=s$, \item[(ii)] $M_i(I:x_n)=M_i(I)$ for all
$i=1,\ldots,s$, and \item[(iii)] $e(S/(I:x_n))=
1/s!\prod_{i=1}^sM_i(I:x_n)$.
\end{enumerate}
Using Noetherian induction and Lemma \ref{purelinear}, condition
(i) and (iii) yield that $S/(I:x_n)$ is Cohen-Macaulay of
codimension $s$ with linear resolution. Since $(I:x_n)$ is
strongly stable we conclude therefore that $(I:x_n)=(x_1,\ldots,
x_s)^a$ for some integer $a$. Now (ii) implies that $M_1(I)=a$, so
that
\[
(x_1,\ldots, x_s)^ax_n\subseteq \mm I\subset I\subset (x_1,\ldots,
x_s)^a,
\]
a contradiction.
\end{proof}

In order to prove the improved multiplicity conjecture in the
quasi-pure case we must assume that the considered algebras are
Cohen-Macaulay. It would be nice this hypothesis could be dropped
from the assumptions.

\begin{Theorem}
\label{almost} Suppose that $I\subset S$ is a Cohen-Macaulay ideal
of codimension $s$ with a quasi-pure resolution. Then $S/I$
satisfies the improved multiplicity conjecture.
\end{Theorem}

\begin{proof} The proof given in \cite{HS} yields the assertion.
We sketch the arguments.   Let
\[ 0\To\bigoplus_{j=1}^{b_s} S(-d_{sj})\To\cdots\To
\bigoplus_{j=1}^{b_1} S(-d_{1j})\To S \To 0.
\]
be the minimal graded free resolution of $S/I$.  There are  square
matrices $A$ and $B$ derived from the shifts $d_{ij}$ with
\begin{eqnarray}
\label{detm} \det A= \sum_{1\leq j_i\leq b_i\atop 1\leq i\leq
s}\prod_{i=1}^sd_{ij_i}V(d_{1j_1},\ldots, d_{sj_s}),
\end{eqnarray}
where $V(d_{1j_1},\ldots, d_{sj_s})$ is the Vandermonde with
entries given by  $d_{1j_1},\ldots, d_{sj_s}$, and such that
\begin{eqnarray}
\label{compare}
 \det A=s!e(S/I)\det B,
 \end{eqnarray}

 It is also shown that
\begin{eqnarray} \label{sum} \det B= \sum_{1\leq j_i\leq b_i\atop
1\leq i\leq s}V(d_{1j_1},\ldots, d_{sj_s}).
\end{eqnarray}
Using the fact that the resolution is quasi-pure, it follows from
(\ref{sum}) that $\det B>0$. Taking minimum and maximum of the
$d_{ij_i}$ in the products in (\ref{detm}) and using
(\ref{compare}), one obtains the inequalities
\begin{eqnarray*}
\label{bounds} (\prod_{i=1}^sm_i)\det B\leq s!e(S/I)\det B
\leq(\prod_{i=1}^sM_i)\det B.
\end{eqnarray*}
Here the lower, resp.\ the upper inequality becomes an  equality
if and only if $m_i=d_{ij}$, resp.\ $M_i=d_{ij}$ for all $i$ and
$j$. Thus the assertion follows.
\end{proof}

As a last example we consider ideals of maximal minors. For these
ideals Mir\'o-Roig has proved Conjecture  1. Inspecting the
inequalities in her proof one can also see that improved
multiplicity conjecture holds. For the  convenience of the reader
we give a complete proof of the theorem, similar to that of
Mir\'o-Roig, in order to see explicitly that the bounds are
reached only when the resolution of the ideal is pure.
Independently, also Migliore, Nagel and R\"omer found a proof of
this theorem.

Let $H=(h_{ij})$ be an $m\times n$-matrix with $m\leq n$, whose
entries are polynomials. We say that $H$ is a {\em homogeneous
matrix} if all minors of $H$ are homogenous polynomials. In
particular, the entries of $H$ itself must be homogenous. For each
$i$ and $j$ let $d_{ij}=\deg h_{ij}$. Then, since the 2-minors are
homogeneous, we get $d_{1j}+d_{i1}=d_{11}+d_{ij}$ for all $i$ and
$j$. Thus if we set $b_i=d_{i1}$ and $a_j=d_{11}-d_{1j}$, then
\[
d_{ij}=b_i-a_j \quad \text{for all}\quad  i=1,\ldots,m\quad
\text{and}\quad j=1,\ldots,n.
\]
Conversely, given any sequences of integers $b_1,\ldots,b_m$ and
$a_1,\ldots,a_n$ we  obtain the degree matrix of a homogeneous
matrix by setting $d_{ij}=b_i-a_j$. After a suitable permutation
of the rows and columns, we may  assume that $a_1\leq a_2\leq
\cdots \leq a_n$ and $b_1\leq b_2  \cdots\leq b_m$. Then this
implies that $d_{i+1,j}\geq d_{ij} \geq  d_{ij+1}$ for all $i$ and
$j$. For the rest of this section we will remain with this
assumption on the degrees of the entries.

Set $r=n-m$, and  let $I_m(H)$ be the ideal of maximal minors of
$M$. Then $\height I_m(H)\leq r+1$, and if equality holds then
$I_m(H)$ is perfect, see \cite[Theorem 2.1 and Theorem 2.7]{BV}.

We want to prove the following

\begin{Theorem}
\label{holds} Suppose that $\height I_m(H)=r+1$. Then the improved
multiplicity conjecture  holds for $S/I_m(H)$.
\end{Theorem}

\begin{proof} For later calculations it is useful to set $u_{ij}=d_{j, i+j-1}$
for all $i=1,\ldots r+1$ and $j=1,\ldots,m$. Using this notation,
we have
\[
b_j-a_{i+j-1}=u_{ij}
\]
for all  $i$ and $j$ in the above range, and since we assume that
$a_1\leq a_2\leq \cdots\leq a_n$ and $b_1\leq b_2\leq \cdots\leq
b_m$ we have
\begin{eqnarray}
\label{inequalities1} u_{1j}\geq u_{2j}\geq \cdots \geq u_{r+1,j}
\quad \text{for}\quad j=1,\ldots, r+1,
\end{eqnarray}
and
\begin{eqnarray}
\label{inequalities2} u_{ij}\geq u_{i+1,j-1}\quad \text{for all
$i$, $j$ with $i+j\leq n+m+1$, $i\leq r+1$ and $1<j$.}
\end{eqnarray}

According to \cite[Corollary 6.5]{HT}  the multiplicity of
$R=S/I_m(H)$ is then given by
\begin{eqnarray}
\label{herzogtrung}
 e(R)=\sum_{1\leq j_1\leq j_2\leq \cdots\leq
j_{r+1}\leq m}\prod_{i=1}^{r+1}u_{i,j_i}.
\end{eqnarray}

Since by assumption  $I_m(H)$ is perfect, the Eagon-Northcott
complex provides a minimal graded free $S$-resolution of $I_m(H)$.
This allows us to compute the numbers $M_i(I_m(H))$.  Following
\cite{BV},  the Eagon-Northcott complex resolving $I_m(H)$ can be
described as follows: let $F$ an $G$ be finitely generated free
$S$-modules with basis $f_1,\ldots, f_m$ and $g_1,\ldots, g_n$,
resp., and let $\phi\:G\to F$ be the linear map with
\[
\phi(g_j)=\sum_{i=1}^mh_{ij}f_i,\quad j=1,\ldots,n.
\]
denote by $\bigwedge^jG$ the $j$th exterior power of $G$, and by
$S_i(F)$ the $i$th symmetric power of $F$. For $j=1,\ldots,n$ we
may then view $\phi(g_j)$ as an element of the symmetric algebra
$S(F)=\Dirsum_iS_i(F)$, and the Koszul complex of
$\phi(g_1),\ldots, \phi(g_n)$ is then given by
\[
0\To \bigwedge^n G\tensor S(F)\To \cdots \To \bigwedge^1 G\tensor
S(F)\To\bigwedge^0 G\tensor S(F)\To 0.
\]
The symmmetric algebra $S(F)$ is graded, and the elements
$\phi(g_j)$ are homogeneous of degree $1$ (the coefficients
$h_{ij}$ are here considered to be of degree $0$). The $r$th
graded component of the Koszul complex

\begin{eqnarray}
\label{koszul} \hspace{0.5cm} 0\To \bigwedge^r G\tensor S_0(F)\To
\cdots \To \bigwedge^1 G\tensor S_{r-1}(F)\To\bigwedge^0 G\tensor
S_r(F)\To 0
\end{eqnarray}
is a complex of free $S$-modules, whose $S$-dual
\begin{eqnarray*}
\label{dualkoszul} (5)\hspace{0.2cm} 0\To (\bigwedge^0 G\tensor
S_r(F))^*\To \cdots \To (\bigwedge^{r-1} G\tensor
S_1(F))^*\To(\bigwedge^r G\tensor S_0(F))^*\To 0
\end{eqnarray*}
is the Eagon-Northcott complex resolving $I_m(H)$.

Set $\deg f_i=a-b_i$ for $i=1,\ldots,m$ and $\deg g_j=a-a_j$ for
$j=1,\ldots,n$. Then the Koszul complex (\ref{koszul}) is a graded
complex. In order to make the augmentation map $(\bigwedge^r
G\tensor S_0(F))^*\to I_m(H)$ homogeneous of degree $0$. We let
the Eagon-Northcott complex be the dual of complex (\ref{koszul})
with respect to $S(-(r-1)a-b)$, where $a=\sum_{i=1}^na_i$ and
$b=\sum_{i=1}^mb_i$.

Now we can compute the $M_i$ for the ideal $I_m(H)$. For a basis
element $e\in \bigwedge^{r-k} G\tensor S_{k}(F)$ we denote by
$e^*$ the  dual basis element  in $(\bigwedge^{r-k} G\tensor
S_{k}(F))^*$. Then the elements
\[
(g_{i_1}\wedge g_{i_2}\wedge \cdots \wedge g_{i_{r-k}}\tensor
f_{j_1}f_{j_2}\cdots f_{j_{k}})^*
\]
establish a basis of $(\bigwedge^{r-k} G\tensor S_k(F))^*$, and we
have
\begin{eqnarray*}
\deg (g_{i_1}\wedge g_{i_2}\wedge \cdots \wedge
g_{i_{r-k}}&\tensor & f_{j_1}f_{j_2}\cdots
f_{j_{k}})^*\\
&=&(r-1)a+b-\sum_{s=1}^{r-k}(a-a_{i_s})-\sum_{t=1}^k(a-b_{j_t})\\
&=&-a+b+\sum_{s=1}^{r-k}a_{i_s}+\sum_{t=1}^kb_{j_t}.
\end{eqnarray*}
It follows that
\begin{eqnarray*}
M_{k+1}&=&-a+b+a_{m+k+1}+\cdots +a_n+kb_m\\
&=&\sum_{i=1}^{m-1}(b_i-a_i)+\sum_{i=0}^k(b_m-a_{m+i})\\
&=&\sum_{j=1}^{m-1}u_{1j}+\sum_{i=1}^{k+1}u_{im}.
\end{eqnarray*}
Thus we need to prove the following inequality
\begin{eqnarray}
\label{upper} (r+1)!\sum_{1\leq j_1\leq j_2\leq \cdots\leq
j_{r+1}\leq m}\prod_{i=1}^{r+1}u_{i,j_i}\leq
\prod_{k=0}^r(\sum_{j=1}^{m-1}u_{1j}+\sum_{i=1}^{k+1}u_{im}).
\end{eqnarray}
We use induction on $\min\{r+1,m\}$ to prove this inequality. In
case $r=0$, we have $n=m$, and on both sides of the inequality we
have  the same expression, namely $\sum_{j=1}^nu_{1j}$. In case
$m=1$, the ideal $I_m(H)$ is generated by the regular sequence
$h_{11},\ldots, h_{1n}$. In this case the inequality is also known
to be true, see \cite{HS}. (It also follows from the result in the
next section).

We now assume that $\min\{r+1,m\}>1$, and decompose the expression
for the multiplicity as follows
\begin{eqnarray*}
(r+1)!\sum_{1\leq j_1\leq j_2\leq \cdots\leq j_{r+1}\leq
m}\prod_{i=1}^{r+1}u_{i,j_i}&=&(r+1)!\sum_{1\leq j_1\leq j_2\leq
\cdots\leq j_{r+1}< m}\prod_{i=1}^{r+1}u_{i,j_i}\\
&+&(r+1)(r!\sum_{1\leq j_1\leq \cdots\leq j_r\leq
m}\prod_{i=1}^ru_{i,j_i})u_{r+1,m}.
\end{eqnarray*}

Using induction we can replace the second summand by the larger
term
\[
(r+1)\prod_{k=0}^{r-1}(\sum_{j=1}^{m-1}u_{1j}+\sum_{i=1}^{k+1}u_{im})u_{r+1,m},
\]
and obtain the inequality
\begin{eqnarray}
\label{right} (r+1)!e(R)&\leq&(r+1)!\sum_{1\leq j_1\leq j_2\leq
\cdots\leq
j_{r+1}< m}\prod_{i=1}^{r+1}u_{i,j_i}\\
&+&(r+1)\prod_{k=0}^{r-1}(\sum_{j=1}^{m-1}u_{1j}+\sum_{i=1}^{k+1}u_{im})u_{r+1,m},\nonumber
\end{eqnarray}
On the other hand, by (\ref{inequalities1}) we get
\begin{eqnarray}\label{left}
\; \; \prod_{k=0}^r(\sum_{j=1}^{m-1}u_{1j}+\sum_{i=1}^{k+1}u_{im})&\geq&\prod_{k=0}^{r-1}(\sum_{j=1}^{m-1}u_{1j}+\sum_{i=1}^{k+1}u_{im})(\sum_{j=1}^{m-1}u_{1j}+(r+1)u_{r+1,m})\\
&=&\prod_{k=0}^{r-1}(\sum_{j=1}^{m-1}u_{1j}+\sum_{i=1}^{k+1}u_{im})(\sum_{j=1}^{m-1}u_{1j})
\nonumber \\
&+&(r+1)\prod_{k=0}^{r-1}(\sum_{j=1}^{m-1}u_{1j}+\sum_{i=1}^{k+1}u_{im})u_{r+1,m}.\nonumber
\end{eqnarray}
Thus comparing (\ref{right}) and (\ref{left}) it remains to be
shown that
\[
(r+1)!\sum_{1\leq j_1\leq j_2\leq \cdots\leq j_{r+1}\leq
m-1}\prod_{i=1}^{r+1}u_{i,j_i}\leq\prod_{k=0}^{r-1}(\sum_{j=1}^{m-1}u_{1j}+\sum_{i=1}^{k+1}u_{im})(\sum_{j=1}^{m-1}u_{1j}).
\]
By induction hypothesis
\[
(r+1)!\sum_{1\leq j_1\leq j_2\leq \cdots\leq j_{r+1}\leq
m-1}\prod_{i=1}^{r+1}u_{i,j_i}\leq\prod_{k=0}^r(\sum_{j=1}^{m-2}u_{1j}+\sum_{i=1}^{k+1}u_{i,m-1}).
\]
Thus the desired inequality follows once we can show that
\begin{eqnarray}\label{end}
\prod_{k=0}^r(\sum_{j=1}^{m-2}u_{1j}+\sum_{i=1}^{k+1}u_{i,m-1})\leq\prod_{k=0}^{r-1}(\sum_{j=1}^{m-1}u_{1j}+\sum_{i=1}^{k+1}u_{im})(\sum_{j=1}^{m-1}u_{1j}).
\end{eqnarray}
This however is obvious, since the $k$th factor  on left hand side
 for $k=0$ is  equal to the last factor on the right hand side, and
since, due to (\ref{inequalities2}), for $k=1,\ldots, r$, the
$k$th factor on the left hand side
$\sum_{j=1}^{m-2}u_{1j}+\sum_{i=1}^{k+1}u_{i,m-1}$ is less than or
equal to the $(k-1)$th factor
$\sum_{j=1}^{m-1}u_{1j}+\sum_{i=1}^{k}u_{im}$ on the left hand
side.

In order to prove the lower inequality, note that
\begin{eqnarray*}
m_{k+1}&=&-a+b+a_{1}+\cdots +a_{r-k}+kb_1\\
&=&\sum_{i=2}^{m}(b_i-a_{r+i})+\sum_{i=0}^k(b_1-a_{r-i+1})\\
&=&\sum_{i=r+1-k}^{r+1}u_{i1}+\sum_{j=2}^{k+1}u_{r+1,j}.
\end{eqnarray*}
Thus we need to prove the following inequality
\begin{eqnarray}
\label{lower}(r+1)!\sum_{1\leq j_1\leq j_2\leq \cdots\leq
j_{r+1}\leq m}\prod_{i=1}^{r+1}u_{i,j_i}\geq
\prod_{k=0}^r(\sum_{i=r+1-k}^{r+1}u_{i1}+\sum_{j=2}^{k+1}u_{r+1,j}).
\end{eqnarray}
The proof of this inequality is completely analogue to that of
inequality (\ref{upper}), since the situation somehow dual to
previous case. Indeed, the substitution $$u_{ij}\mapsto u_{r+2-i,
m+1-j}\quad \text{for}\quad i=1,\ldots,r\quad\text{and}\quad
j=1,\ldots,m$$ transfers (\ref{upper}) to (\ref{lower}) and {\em
reverses} the inequalities (\ref{inequalities1}) and
(\ref{inequalities2}). Thus the lower bound follows from the upper
bound.

Now suppose the multiplicity of $S/I_m(H)$ reaches the upper
bound. This is only possible if we have equality in (\ref{right}),
(\ref{left}) and (\ref{end}).

It follows from the formula for the shifts in the resolution of
the Eagon-Northcott complex, that the resolution of $S/I_m(H)$ is
pure if and only if $a_1=a_2=\cdots =a_n$ and $b_1=b_2=\cdots
=b_m$, which is equivalent to say that all $u_{ij}$ are equal.

Hence by induction we have equality in (\ref{right}) if and only
if $u_{ij}=u$ for some $u$ and all $i=1,\ldots, r+1$ and
$j=1,\ldots, m-1$. On the other hand, we get equality in
(\ref{left}) if and only if $u_{im}=u_{r+1,m}$ for $i=1,\ldots,r$,
while equality holds in (\ref{end}) if and only if
$u_{im}=u_{i,m-1}=u$ for $i=1,\ldots,r$. Thus all $u_{ij}$ must be
equal to $u$.

Using the reflection principle of above it also follows that the
lower bound for the multiplicity is reached only when the
resolution of $S/I_m(H)$ is pure.
\end{proof}

\begin{Remarks}
{\em  (a) Theorem \ref{holds} includes the case studied by Guardo
and Van Tuly \cite{GT}, namely that rings whose defining ideal is
generated by powers of a homogeneous regular sequence satisfy the
improved multiplicity conjecture. In fact, if $f_1,\ldots,
f_{r+1}$ is a homogeneous regular sequence and $I=(f_1,\ldots,
f_{r+1})$, then $I^m$ is the ideal of maximal minors of the
$m\times m+r$-matrix whose $i$th main diagonal has all entries
$f_i$ for $i=1,\ldots, r+1$, while all other entries of the matrix
are $0$.

(b) It can be easily seen from the proof of Theorem \ref{regular}
that $e(R/({\bf f}))$  reaches the upper bound if and only if
there exists an integer $d$ such that $M_i=id$ and $\deg f_i=d$
for all $i$. }
\end{Remarks}


\newpage
\begin{thebibliography}{99}
\bibitem{AHH} A.\ Aramova, J.\ Herzog and T.\ Hibi, Ideals with
stable Betti numbers, Adv.\ Math. {\bf 152}(1), (2000), 72--77.


\bibitem{AHH1} A.\ Aramova, J.\ Herzog and T.\ Hibi, Shifting
operations and graded Betti-numbers, J. Algebraic Combin.\ {\bf
12} (2000), 207--222.

\bibitem{BS} D.\ Bayer and M.\ Stillman, A criterion or detecting
$m$-regularity, Invent.\ Math.\ {\bf 87}(1), (1987), 1--11.




\bibitem{BH} W.\ Bruns and J.\ Herzog, {\em Cohen--Macaulay
rings},  Revised Edition, Cambridge University Press, Cambridge,
1996.


\bibitem{BV} W.\ Bruns and U.\ Vetter, {\em Determinantal rings},
Springer Verlag, Graduate texts in Mathematics {\bf 150} (1995).


\bibitem{CHT} D.\ Cutkosky, J.\ Herzog and N.V.\ Trung, Asymptotic
behaviour of the Castelnuovo-Mumford regularity, Compositio Math.
{\bf 118} (1999), no. 3, 243--261.



\bibitem{ERT}  D.\ Eisenbud, A.\ Reeves and B.\ Totaro, Initial
ideals, Veronese subrings, and rates of Algebras, Adv.\ Math.\
{\bf 109} (1994), 168--187.

\bibitem{EK} S.\ Eliahou and M.\ Kervaire, Minimal resolutions for
some monomial ideals, J.\ Algebra {\bf 129}, (1990), 11--25.

\bibitem{GT} E.\ Guardo and A.\ Van Tuyl, Powers of complete
intersections: graded Betti-numbers and applications, Preprint
2004. math.AC/0409090.

\bibitem{G} L.H.\ Gold, A degree bound for codimension two lattice
ideals, J.\ Pure Appl.\ Algebra {\bf 182} (2003) 201--207.



\bibitem{HH} J.\ Herzog and T.\ Hibi, Componentwise linear ideals,
Nagoya Math.\ J.\ {\bf 153}, (1999), 141--153.

\bibitem{HK} J.\ Herzog and M.\ K\"uhl, On the Betti-numbers of
finite pure and linear resolutions. Comm.\ Alg.\ {\bf 12}  (1984),
1627--1646.



\bibitem{HS} J.\ Herzog and H.\ Srinivasan, Bounds for
multiplicities, Trans. AMS {\bf 350} (1998), 2879--2902.


\bibitem{HT} J.\ Herzog and N.V.\ Trung, Gr\"obner bases and
multiplicity of determinantal and pfaffian ideals, Adv.\ Math.
{\bf 96}, (1992), 1--37.



\bibitem{HM} C.\ Huneke and M.\ Miller, A note on the multiplicity
of Cohen-Macaulay algebras with pure resolution, Canad.\ J.\
Math.\ {\bf 37} (1985), 1149--1162.


\bibitem{K}  V.\ Kodiyalam, Asymptotic behaviour of Castelnuovo-Mumford
regularity, Proc. Amer. Math. Soc. {\bf 128} (2000), no. 2,
407--411


\bibitem{MNR} J.\ Migliore, U.\ Nagel and T.\ R\"omer, The
multiplicity conjecture in low codimension, Preprint 2004.
math.AC/0410497

\bibitem{M} R.M.\ Mir\'o-Roig, A note on the multiplicity conjecture of
determinantal rings, Preprint 2005. math.AC/0504077


\bibitem{R} T.\ R\"omer, Note on bounds for multiplicities, J.\
Pure  Appl.\ Algebra {bf 195} (2005), 113--123.


\bibitem{S} H.\ Srinivasan, A note on the multiplicities of
Gorensetin algebras, J.\ Algebra {\bf 208}(2) (1998), 425--443.
\end{thebibliography}
\end{document}